 \documentclass[a4paper,11pt]{article}
  \usepackage{amssymb,amsbsy,amsmath,amscd, latexsym,theorem,eepic,enumerate}
  \usepackage{xypic,times,color}
  \xyoption{curve}
  \usepackage{geometry}
{\theorembodyfont{\rmfamily}\newtheorem{example}{Example}[section]}
{\theorembodyfont{\rmfamily}\newtheorem{Def}[example]{Definition}}
{\theorembodyfont{\rmfamily}\newtheorem{cond}[example]{Conditions}}

\newtheorem{Prop}[example]{Proposition}

{\theorembodyfont{\rmfamily}}
\newtheorem{corollary}[example]{Corollary}

\newtheorem{blank}[example]{\hspace{-0.45em}}

\newtheorem{question}[example]{Question}

\newenvironment{pf}{\noindent {\bf Proof}}{\rule{0mm}{1mm}\hfill $\Box$

\mbox{}}
\newenvironment{enumeratein}{\begin{enumerate}[\vspace{-1ex}\hspace{1em}
(a)]}{\end{enumerate}}

\def\Ob{\mathrm{Ob}}

\def\Set{\mathsf{Set}}
\def\Fun{\mathsf{FUN}}
\def\cls{\mathsf{cls}\,}

\def\orbitgpd{/\!/}

\def\~{\sim}
\def\ab{\mathrm{ab}}

\newcommand{\labto}[1]{\stackrel{#1}{\longrightarrow}}
\def\I{\mathbf{I}}
\def\ES{\mathbf{S}}
\def\<{\langle}
\def\>{\rangle}
\def\#{\sharp}

\def\midsq{\ar @{} [dr] |}

\def\St{\mathrm{St}}
\newcommand{\Real}{{\mathbb R}}

\newcommand{\threeaxes}[3]{\def\objectstyle{\scriptstyle}  \objectmargin={0pt}
\xy
(0,0)*+{}="a",(0,-6)*+{\rule{0em}{1.5ex}#2}="b",(7,0)*+{\;#1}="c",
(14,-3)*+{\;#3}="d" \ar@{->} "a";"b" \ar @{->}"a";"c"  \ar
@{->}"a";"d"\endxy }

\newcommand{\directs}[2]{\def\objectstyle{\scriptstyle}  \objectmargin={0pt}
\xy
(0,4)*+{}="a",(0,-2)*+{\rule{0em}{1.5ex}#2}="b",(7,4)*+{\;#1}="c"
\ar@{->} "a";"b" \ar @{->}"a";"c" \endxy }
\newcommand{\xdirects}[2]{\def\objectstyle{\scriptstyle}  \objectmargin={0pt}
\xy
(0,0)*+{}="a",(0,-6)*+{\rule{0em}{1.5ex}#2}="b",(7,0)*+{\;#1}="c"
\ar@{->} "a";"b" \ar @{->}"a";"c" \endxy }
\newcommand{\sdirects}[2]{\def\objectstyle{\scriptstyle}  \objectmargin={0pt}
\xy
(0,2.2)*+{}="a",(0,-2.5)*+{\rule{0em}{1.5ex}#2}="b",(7,2.2)*+{\;#1}="c"
\ar@{->} "a";"b" \ar @{->}"a";"c" \endxy }

\def\rho{\varrho}

\def\epsilon{\varepsilon}

\def\Z{\mathbb{Z}}

\def\lan{\langle}
\def\ran{\rangle}

\def\epsilon{\varepsilon}

\def\Z{\mathbb{Z}}

\def\Ker{\mathop{\rm Ker}\nolimits}

\def\epsilon{\varepsilon}

\newcommand{\io}{^{-1}}
\begin{document}

\title{The fundamental groupoid of the quotient  of a Hausdorff
space\\ by a discontinuous action of a discrete group\\ is the
orbit groupoid of the induced action}
\author{Ronald Brown\thanks{email: r.brown@bangor.ac.uk}, \\
 Mathematics Division, \\ School of Informatics,  \\ University of
Wales, Bangor \\Gwynedd LL57 1UT,  U.K. \and  Philip
J. Higgins\thanks{email: p.j.higgins@durham.ac.uk}, \\
Department of Mathematical Sciences, \\ Science Laboratories, \\ South Rd., \\
Durham, DH1 3LE,  U.K.} \maketitle

\begin{center}
University of Wales, Bangor, Maths Preprint  02.25
\end{center}
\begin{abstract}
The main result is that the fundamental groupoid of the orbit
space of a discontinuous  action of a discrete group on a
Hausdorff space which admits a universal cover is the orbit
groupoid of the fundamental groupoid of the space. We also
describe work of Higgins and of Taylor which makes this result
usable for calculations. As an example, we compute the fundamental
group of the symmetric square of a space.

The main result, which is related to work of Armstrong, is due to
Brown and Higgins in 1985 and was published in sections 9 and 10
of Chapter 9 of the first author's book on Topology
\cite{Brown:1988}. This is a somewhat edited, and in one point (on
normal closures) corrected, version of those sections. Since the
book is out of print, and the result seems not well known, we now
advertise it here.

It is hoped that this account  will also allow wider views of
these  results, for example in topos theory and descent theory.

Because of its provenance, this should be read as a graduate text
rather than an article. The Exercises should be regarded as
further propositions for which we leave the proofs to the reader.
It is expected that this material will be part of a new edition of
the book.
\end{abstract}
{\small  \noindent MATH CLASSIFICATION: 20F34, 20L13, 20L15,
57S30}

\section{Groups acting on spaces}
In this section we show some of the theory of a group $G$ acting
on a topological space $X$, and describe the {\em orbit
topological space,} which is written $X/G$.

There arises the problem of relating topological invariants of the
orbit  space $X/G$ to those of $X$ and the group action.  In
particular, it is a complicated and interesting question to find,
 if at all possible, relations between the fundamental groups and
groupoids of $X$ and $X/G$.  This we shall do for a particular
family of actions which arise commonly, namely the
{\em discontinuous actions.}  The resulting % page 9.9-2
theory generalises that of regular covering spaces, and has a
number of important applications.  A useful special case of a
discontinuous action is the action of a finite group on a
Hausdorff space (see below); there are in the literature many
interesting cases of discontinuous actions of infinite groups (see
\cite{Beardon1}).

We now come to formal definitions.

Let $G$ be a group, with its group operation written as
multiplication, and let $X$ be a set.  An {\em action} of $G$ on
$X$ is a function $G \times X \to X $, written $(g,x) \mapsto g
\cdot x$, satisfying the following properties for all $x$ in $X$
and $g,h$ in $G$:

\begin{blank}  \label{P:9.9.1} {\em (i)} $1 \cdot x = x, $\\
{\rule{1.47em}{0ex}} {\em (ii) } $g\cdot(h \cdot x) = (gh) \cdot
x. $
\end{blank}

Thus the first rule says that the identity of $G$ acts as
identity, and the second rule says that two elements of $G$,
acting successively, act as the product of the two elements.

There are some standard notions associated with such an action.
First, an equivalence relation is defined on $X$ by $x\sim y$ if
and only if there is an element $g$ of $G$ such that $y = g \cdot
x$. This is an equivalence relation.  Reflexivity follows since
$G$ has an identity.  Symmetry follows from the existence of
inverses in $G$, using \ref{P:9.9.1} (i), (ii). Transitivity
follows from the product of two elements in $G$ being in $G$.  The
equivalence classes under this relation are the {\em orbits} of
the action.  The set of these orbits is written $X/G$.

Suppose given an action of the group $G$ on the set $X$.  If $x
\in X$, then the {\em group of stability of} $x$ is the subgroup
of $G$ \[ G_{x} = \{g \in  G : g \cdot x = x\}  \] The elements of
$G_{x}$ are said to {\em stabilise $x$,} that is, they leave $x$
fixed by their action.  If $G_{x}$ is the whole of $G$, then $x$
is said to be a {\em fixed point} of the action.  The set of fixed
points of the action is often written $X^{G}$.  The action is said
to be {\em free} if all groups of stability are trivial.

Another useful condition is the notion of {\em effective} action
of a group.  This requires that if $g,h \in G$ and for all $x$ in
$X, \, g \cdot x = h \cdot x $, then $g = h$.  In this case the
elements of $G$ are entirely determined by their action on $X$.

We now turn to the topological situation.  Let $X$ be
% page 9.9-2
a topological space, and let $G$ be a group.  An {\em action} of
$G$ on $X$ is again a function $G \times X \to X$ with the same
properties as given in \ref{P:9.9.1}, but with the additional
condition that when $G$ is given the discrete topology then the
function $(g,x) \mapsto g \cdot x$ is continuous.  This amounts to
the same as saying that for all $g \in G$, the function $g_{\#} :
x \mapsto g \cdot x$ is continuous. Note that $g_{\#}$ is a
bijection with inverse $(g^{-1})_{\#}$, and since these two
functions are continuous, each is a homeomorphism.

Let $X/G$ be the set of orbits of the action and let $p : X \to
X/G$ be the {\em quotient map}, which assigns to each $x$ in $X$
its orbit. For convenience we will write the orbit of $x$ under
the action as $\bar{x}$.  So the defining property is that
$\bar{x} = \bar{y}$ if and only if there is a $g$ in $G$ such that
$y = g \cdot x$.  Now a topology  has  been given for  $X$.  We
therefore give the orbit space $X/G$ the identification topology
with respect to the map $p$.  This topology will always be assumed
in what follows.  The first result on this topology, and one which
is used a lot, is as follows.

\begin{Prop} \label{P:9.9.3} The quotient map $p : X \to X/G$ is
an open map.   \end{Prop}

\begin{pf} Let $U$ be an open set of $X$.  For each $g \in G$
the set $g\cdot U$, by which is meant the set of $g \cdot x$ for
all $x$ in $U$, is also an open set of $X$, since $g_{\#}$ is a
homeomorphism, and $g\cdot U = g_{\#}[U]$.  But \[ p^{-1}p[U] =
\bigcup_{g \in G}g\cdot U .\] Since the union of open sets is
open, it follows that $p^{-1}p[U]$ is open, and hence $p[U]$ is
open.
\end{pf}
\begin{Def}
An action of the group $G$ on the space $X$ is called {\em
discontinuous} if the stabiliser of each point of $X$ is finite,
and each point $x$ in $X$ has a neighbourhood $V_{x}$ such that
any element $g$ of $G$ not in the stabiliser of $x$ satisfies
$V_{x} \cap g\cdot V_{x} = \emptyset$. \end{Def} Suppose $G$ acts
discontinuously on the space $X$.  For each $x$ in $X$ choose such
an open neighbourhood $V_{x}$ of $x$.  Since the stabiliser
$G_{x}$ of $x$ is finite, the set \[ U_{x} = \bigcap \{g\cdot
V_{x} : g \in G_{x}\} \] is open; it contains $x$ since the
elements of $G_{x}$ stabilise $x$. Also if $g \in G_{x}$ then
$g\cdot U_{x} = U_{x}$.  We say $U_{x}$ is {\em invariant} under
the action of the group $G_{x}$.  On the other hand, if $h \not\in
G_{x}$ then
$$(h\cdot U_{x}) \cap U_{x} \subseteq (h\cdot V_{x}) \cap V_{x} =
\emptyset .$$

An open neighbourhood $U$ of $x$ which satisfies $(h\cdot U) \cap
U=  \emptyset $ for $h \not\in G_{x}$ and is invariant under the
action of $G_{x}$ is called a {\em
canonical % page 9.9-4
neighbourhood} of $x$.  Note that any neighbourhood $N$ of $x$
contains a canonical neighbourhood: the proof is obtained by
replacing $V_{x}$ in the above by $N \cap V_{x}$.  The image in
$X/G$ of a canonical neighbourhood $U$ of $x$ is written
$\overline{U}$ and called a {\em canonical neighbourhood} of
$\bar{x}$.

In order to have available our main example of a discontinuous
action, we prove:

\begin{Prop} \label{P:9.9.4} An action of a finite group on a
Hausdorff space is discontinuous. \end{Prop}

\begin{pf} Let $G$ be a finite group acting on the Hausdorff
space $X$.  Then the stabiliser of each point of $X$ is a subgroup
of $G$ and so is finite.

Let $x \in X$.  Let $x_{0}, x_{1}, \dots, x_{n}$ be the distinct
points of the orbit of $x$, with $x_{0} = x$.  Suppose $x_{i} =
g_{i} \cdot x, g_{i} \in G, i = 1,\dots, n$, and set $g_{0} = 1$.
Since $X$ is Hausdorff, we can find pairwise disjoint open
neighbourhoods $N_{i}$ of $x_{i}, i = 0,\dots, n$.  Let \[ N =
\bigcap \{g^{-1}_{i}\cdot N_{i} : i = 0,1,\dots n\} \] Then $N$ is
an open neighbourhood of $x$.  Also, if $g \in G$ does not belong
to the stabiliser $G_{x}$ of $x$, then for some $j = 1,\dots, n, g
\cdot x = x_{j}$, whence $g\cdot N \subseteq N_{j}$.  Hence $N
\cap  g\cdot N = \emptyset$, and the action is discontinuous.
\end{pf}  Our main result in general topology on discontinuous
actions is the following.

\begin{Prop} \label{P:9.9.5} If the group $G$ acts
discontinuously on the Hausdorff space $X$, then the quotient map
$p : X \to X/G$ has the path lifting property:  that is, if
$\bar{a} : \I \to X/G$ is a path in $X/G$ and $x_{0}$ is a point
of $X$ such that $p(x_{0}) = \bar{a}(0)$, then there is a path $a
: \I \to X$  such that $pa = \bar{a}$  and $a(0) = x_{0}$.
\end{Prop}

\begin{pf} If there is a lift $a$ of $\bar{a}$ then there is an
element $g$ of $G$ such that $g\cdot a(0) = x_{0}$ and so $g\cdot
a$ is a lift of $\bar{a}$ starting at $x_{0}$.  So we may ignore
$x_{0}$ in what follows.

Since the action is discontinuous, each point $\bar{x}$ of $X/G$
has a canonical neighbourhood.  By the Lebesgue covering lemma,
there is a subdivision
\[ \bar{a} = \bar{a}_{n} + \dots  + \bar{a}_{1} \] of $\bar{a}$
such that the image of each $\bar{a}_{i}$ is contained in a
canonical neighbourhood.  So if the path lifting property holds
for each canonical neighbourhood in $X/G$, then it holds for
$X/G$.

Now $p : X \to X/G$ is an open map.  Hence for all $x$ in $X$
the % page 9.9-5
restriction $p_{x} : U_{x} \to \overline{U}_{x}$ of $p$ to a
canonical neighbourhood $U_{x}$ is also open, and hence is an
identification map.  So we can identify $\overline{U}_{x}$ with
the orbit space $(U_{x})/G_{x}$.  The key point in this case is
that the group $G_{x}$ is finite.

Thus it is sufficient to prove the path lifting property for the
case of the action of a {\em finite} group $G$, and this we do by
induction on the order of $G$.  That is, we assume that path
lifting holds for any action of any proper subgroup of $G$ on a
Hausdorff space, and we prove that path lifting holds for the
action of $G$. The case $|G|=1$ is trivial.

Again let $\bar{a}$ be as in the proposition, and we are assuming
$G$ is finite.  Let $F$ be the set of fixed points of the action.
Then $F$ is the intersection for all $g \in G$ of the sets $X^{g}
= \{x \in X : g \cdot x = x\}$. Since $X$ is Hausdorff, the set
$X^{g}$ is closed in $X$, and hence $F$ is closed in $X$.  So
$p[F]$ is closed, since $p^{-1}p[F] = F$.   Let $A$ be the
subspace of $\I$ of points $t$ such that $\bar{a}(t)$ belongs to
$p[F]$, that is, $A = \bar{a}^{-1}p[F]$. Then $A$ is closed.

The restriction of the quotient map $p$ to $p' : F \to p[F]$ is a
homeomorphism.  So $\bar{a}\mid_{A}$ has a unique lift to a map
$a\mid_{A} : $A$ \to X$.  So we have to show how to lift
$\bar{a}\mid_{(\I  \backslash A)}$ to give a map $a\mid_{\I
\backslash A}$ and then show that the function $a : \I \to X$
defined by these two parts  is continuous.

In order to construct $a\mid _{\I \backslash A}$, we first assume
$A = \{1\}$.

Let $S$ be the set of $s \in \I$ such that $\bar{a}\mid_{[0,s]}$
has a lift to a map $a_{s}$ starting at $x$.  Then $S$ is
non-empty, since $0 \in S$.  Also $S$ is an interval.  Let $u =
\sup S$.  Suppose $u<1$.  Then there is a $y \in X \backslash F$
such that $p(y) = \bar{a}(u)$.  Choose a canonical neighbourhood
$U$ of $y$.  If $u>0$,  there is a $\delta > 0$ such that
$\bar{a}[u-\delta, u+\delta ] \subseteq \bar{U}$.  Then $u-\delta
\in S$ and so there is a lift $a_{u-\delta}$ on $[0,u-\delta ]$.
Also the stabiliser of $y$ is a proper subgroup of $G$, since $y
\notin F$, and so by the inductive assumption there is a lift of
$\bar{a}\mid_{[u-\delta, u+\delta]}$ to a path starting at
$a_{u-\delta}(u-\delta)$. Hence we obtain a lift $a_{u+\delta}$,
contradicting the definition of $u$. We get a similar
contradiction to the case $u=0$ by replacing in the above
$u-\delta$ by $0$.  It follows that $u = 1$. We are not quite
finished because all we have thus ensured is that there is a lift
on $[0,s]$ for each $0\leqslant s<1$, but this is not the same as
saying that there is a lift on $[0,1)$. We now prove that such a
lift exists.

By the definition of $u$, and since $u = 1$, for each integer
$n\geqslant 1$ there is a lift $a^{n}$ of
$\bar{a}\mid_{[0,1-n^{-1}]}$. Also there is an element $g_{n}$
of $G$ such that % page 9.9-6
\[ g_{n}\cdot a^{n+1}(1-n^{-1}) =
a^{n}(1-n^{-1})  \] Hence $a^{n}$ and
$g_{n}\cdot(a^{n+1}\mid_{[1-n^{-1},1-(n+1)^{-1}]})$ define a lift
of $\bar{a}\mid_{[0,1-(n+1)^{-1}]}$.  Starting with $n = 1$, and
continuing in this way, gives a lift of $\bar{a}\mid_{[0,1)}$.
This completes the construction of the lift on $\I \backslash A$
in the case $A = \{1\}$.

We now construct a lift of $\bar{a}\mid_{\I \backslash A}$ in the
general case.  Since $A$ is closed, $\I \backslash A$ is a union
of disjoint open intervals each with end points in $\{0,1\} \cup
A$. So the construction of the lift is obtained by starting at the
mid point of any such interval and working backwards and forwards,
using the case $A = \{1\}$, which we have already proved.

The given lift of $\bar{a}\mid _A$ and the choice of lift of
$\bar{a}\mid_{\I \backslash A}$ together define a lift $a : \I \to
X$ of $\bar{a}$ and it remains to prove that $a$ is continuous.

Let $t \in \I$.  If $t \notin A$, then $a$ is continuous at $t$ by
construction.  Suppose then $t \in A$ so that $y = a(t) \in F$.
Let $N$ be any neighbourhood of $y$.  Then $N$ contains a
canonical neighbourhood $U$ of $y$.  If $g \in G$ then $g \cdot y
= y$, and so $U$ is invariant under the action of $G$.  Hence
$p^{-1}p[U] = U$. Since $\bar{a}$ is continuous, there is a
neighbourhood $M$ of $t$ such that $\bar{a}[M] \subseteq \bar{U}$.
Since $pa = \bar{a}$ it follows that $a[M] \subseteq
p^{-1}[\bar{U}] = U$.  This proves continuity of $A$, and the
proof of the proposition is complete. \end{pf}

In our subsequent results we shall use the path lifting property
rather than the condition of the action being discontinuous.

Our aim now is to determine the fundamental groupoid of the orbit
space $X/G$.  In general it is difficult to say much. However we
can give reasonable and useful conditions for which the question
can be completely answered.  From our point of view the result is
also of interest in that our statement and proof use groupoids in
a crucial way.  This use could be overcome, but at the cost of
complicating both the statement of the theorem and its proof.

In order to make the transition from the topology to the algebra,
it is necessary to introduce the notion of a group acting on a
groupoid.

Let $G$ be a group and let $\Gamma$ be a groupoid.  We will write
the group structure on $G$ as multiplication, and the
groupoid % page 9.9-7
structure on $\Gamma$ as addition.  An {\em action} of $G$ on
$\Gamma$ assigns to each $g \in G$ a morphism of groupoids $g_{\#}
: \Gamma \to \Gamma$ with the properties that $1_{\#} = 1 : \Gamma
\to \Gamma$, and if $g,h \in G$ then $(hg)_{\#} = h_{\#}g_{\#}$.
If $g \in G, x \in \Ob(\Gamma), a \in \Gamma$, then we write $g
\cdot x$ for $g_{\#}(x), g\cdot a$ for $g_{\#}(a)$.  Thus the
rules \ref{P:9.9.1}  apply also to this situation, as well as the
laws \noindent {\bf 1.1} (iii) $g.(a+b)=g.a + g.b$, and (iv)
$g.0_x=0_{g.x}$ \\ for all $g \in G, x \in \Ob(\Gamma) , a,b \in
Gamma $ such that $a+b$ is defined.

The action of $G$ on $\Gamma$ is {\em trivial} if $g_{\#} = 1$ for
all $g$ in $G$.

% NB:  The \orbitgpd really ought to be slanted to the right,
%but % I can't find a AMSLaTeX character like that. [I have
%made a command which gives this!]

\begin{Def} Let $G$ be a group acting on a groupoid $\Gamma$.  An {\em orbit
groupoid} of the action is a groupoid $\Gamma \orbitgpd G$
together with a morphism $p : \Gamma \to \Gamma \orbitgpd G$ such
that: \label{P:9.9.6} \begin{enumerate}[(i)] \item
\label{P:9.9.6:i} If $g \in  G, \gamma  \in \Gamma$, then
$p(g\cdot\gamma) = p(\gamma)$. \item \label{P:9.9.6:ii}  The
morphism $p$ is universal for (i), {\it i.e.}  if $\phi: \Gamma
\to \Phi$ is a morphism of groupoids such that $\phi
(g\cdot\gamma) = \phi (\gamma)$ for all $g \in G, \gamma \in
\Gamma$, then there is a unique morphism $\phi^{*} : \Gamma
\orbitgpd G \to \Phi $ of groupoids such that $\phi^{*}p = \phi$.
\end{enumerate}

The morphism $p : \Gamma \to \Gamma \orbitgpd G$ is then called an
{\em orbit morphism}. \hfill $\Box$\end{Def}

The universal property (ii) implies that $\Gamma \orbitgpd G$, if
it exists, is unique up to a canonical isomorphism. At the moment
we are not greatly concerned with proving any general statement
about the {\em existence} of the orbit groupoid. One can argue
that $\Gamma \orbitgpd G$ is obtained from $\Gamma$ by imposing
the relations $g\cdot\gamma = \gamma$ for all $g \in G$ and all
$\gamma \in \Gamma$; however we have not yet explained quotients
in this generality.  We will later prove existence by giving  a
{\em construction} of $\Gamma \orbitgpd G$ which will be useful in
interpreting our main theorem.  But our next result will give
conditions which ensure that the induced morphism $\pi X \to \pi
(X/G)$ is an orbit morphism, and our proof will not {\em assume}
general results on the existence of the orbit groupoid. The reason
we can do this is that our proof directly verifies a universal
property.

First we must point out that if the group $G$ acts on the space
$X$, then $G$ acts on the fundamental groupoid $\pi X$, since each
$g$ in $G$ acts as a homeomorphism of $X$ and $g_{\#} : \pi X \to
\pi X$ may be defined to be the induced morphism.  This is one
important advantage of groupoids over groups: by contrast, the
group $G$ acts on the fundamental group $\pi (X,x)$ only if $x$ is
a fixed point of the action.

Suppose now that $G$ acts on the space $X$.  Our purpose is to
give conditions on the action which enable us to prove that \[
p_{*} : \pi X \to \pi (X/G) \] determines an isomorphism $(\pi
X)\orbitgpd G \to \pi (X/G)$, by verifying the universal property
for $p_{*}$.  We require the following conditions:

\begin{cond} \label{P:9.9.7} \begin{enumerate}[(i)] \item
\label{P:9.9.7:i} The projection $p : X \to X/G$ has the path
lifting property: {\it i.e.}  if $\bar{a} : I \to X/G$ is a path,
then there is a path $a : I \to X$ such that $pa = \bar{a}$.
% page 9.9-8
\item \label{P:9.9.7:ii} If $x \in  X$, then $x$ has an open
neighbourhood $U_{x}$ such that
\begin{enumerate}[(a)]
    \item \label{P:9.9.7:ii:a} if $g \in G$ does not belong to the
stabiliser $G_{x}$ of $X$, then $U_{x} \cap (g\cdot U_{x}) =
\emptyset$;
    \item \label{P:9.9.7:ii:b} if $a$ and $b$ are paths in $U_{x}$
beginning at $x$ and such that $pa$ and $pb$ are homotopic rel end
points in $X/G$, then there is an element $g \in G_{x}$ such that
$g\cdot a$ and $b$ are homotopic in $X$ rel end points.
\end{enumerate} \end{enumerate} \end{cond}
$$\def\labelstyle{\textstyle} 
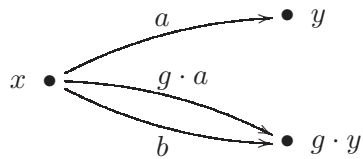
\xymatrix@R=1pc @C=5pc{ & **[r] \bullet \;\; y \\
**[l] x \;\; \bullet \ar @/^/ [ur] ^a \ar @/^/ [dr] ^{g \cdot a}
\ar @/_/ [dr] _b &
\\ & **[r] \bullet \;\; g\cdot y }$$

For a discontinuous action, \ref{P:9.9.7}(ii:a) trivially holds,
while \ref{P:9.9.7}(i) holds by virtue of \ref{P:9.9.5}. However,
\ref{P:9.9.7}(ii:b) is an extra condition.  It does hold if $X$ is
semi-locally simply-connected, since then for sufficiently small
$U$ and $x,g \cdot y \in U$, any two paths in $U$ from $x$ to $g
\cdot y $ are homotopic in $X$ rel end points; so
\ref{P:9.9.7}(ii:b) is a reasonable condition to use in connection
with covering space theory.

A neighbourhood $U_{x}$ of $x \in X$ satisfying \ref{P:9.9.7}(i)
and \ref{P:9.9.7}(ii) will be called a {\em strong canonical
neighbourhood} of $X$.  The image $p[U_{x}]$ of $U_{x}$ in $X/G$
will be called a {\em strong canonical neighbourhood} of $px$.

\begin{Prop} \label{P:9.9.8} If the action of $G$ on $X$
satisfies \ref{P:9.9.7}(i) and \ref{P:9.9.7}(ii) above, then the
induced morphism $p_{*} : \pi X \to \pi (X/G)$  makes $\pi (X/G)$
the orbit groupoid of $\pi X $ by the action of $G$.
\end{Prop}

\begin{pf} Let $\phi : \pi X \to \Phi $ be a morphism to a
groupoid $\Phi$ such that $\phi (g\cdot\gamma) = \phi (\gamma)$
for all $\gamma \in \pi X$ and $g \in G$.  We wish to construct a
morphism $\phi^{*} : \pi (X/G) \to \Phi $ such that $\phi^{*}p =
\phi$.

Let $\bar{a}$ be a path in $X/G$.  Then $\bar{a}$ lifts to a path
$a$ in $X$. Let $[b]$ denote the homotopy class rel end points of
a path $b$.  We prove that $\phi [a]$ in $\Phi$ is independent of
the choice of $\bar{a}$ in its homotopy class and of the choice of
lift $a$; hence we can define $\phi^{*}[\bar{a}]$ to be $\phi
[a]$.

Suppose given two homotopic paths $\bar{a}$ and $\bar{b}$ in
$X/G$, with lifts $a$ and $b$ which without loss of generality we
may assume start at the same point $x$ in $X$.  (If they do not
start at the same point, then one of them may be translated by the
action of $G$ to start at the same point as the other.) Let $h :
\I \times \I \to X/G$ be a homotopy rel end points $\bar{a} \simeq
\bar{b}$.  The method now is not to lift the
homotopy % page 9.9-9
$h$ itself, but to lift pieces of a subdivision of $h$; it is here
that the method differs from that used in the theory of covering
spaces given in Section 9.1 of \cite{Brown:1988}.

Subdivide $\I \times \I$, by lines parallel to the axes, into
small squares each of which is mapped by $h$ into a strong
canonical neighbourhood in $X/G$.  This subdivision determines a
sequence of homotopies $h_{i} : \bar{a}_{i-1} \simeq \bar{a}_{i},
i = 1,2,\dots,n$, say, where $\bar{a}_{0} = \bar{a}, \bar{a}_{n} =
\bar{b}$.  Keep $i$ fixed for the present.  Each $h_{i}$ is
further expressed by the subdivision as a composite of homotopies
$h_{ij}(j = 1,2,\dots,m)$ as shown in the following picture in
which for convenience the boundaries of the $h_{ij}$ are labelled:
$$  \spreaddiagramrows{3.0pc}
\spreaddiagramcolumns{3.0pc}\objectmargin{0pc} \objectwidth{0pc}
\def\labelstyle{\textstyle}
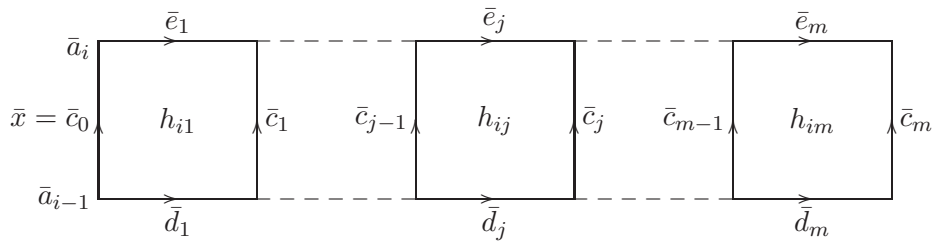
\diagram   \rline^{\bar{e}_1}|\tip
\midsq{h_{i1}} &  \rdashed   &  \midsq{h_{ij}}
\rline^{\bar{e}_j}|\tip& \rdashed  &  \rline^{\bar{e}_m}|\tip
\midsq{h_{im}} & \\
\uline^(0){\bar{a}_{i-1}}^{\bar{x}=\bar{c}_0}|\tip^>>{\bar{a}_{i}
 } \rline_{\bar{d}_1}|\tip&  \rdashed  \uline_{\bar{c}_1}|\tip
&\uline^{\bar{c}_{j-1}}|\tip \rline_{\bar{d}_j}|\tip
&\uline_{\bar{c}_j}|\tip  \rdashed & \uline|\tip^{\bar{c}_{m-1}}
\rline_{\bar{d}_m}|\tip &  \uline_{\bar{c}_m}|\tip \enddiagram $$

Choose lifts $a_{i-1}, a_{i}$ of $\bar{a}_{i-1}, \bar{a}_{i}$
respectively; express $a_{i-1}$ as a sum $a_{i-1} = d_{m} + \dots
+ d_{1}$ and $a_{i}$ as a sum $a_{i} = e_{m} + \dots + e_{1}$
where $d_{j}$ lifts $\bar{d}_{j}$ and $e_{j}$ lifts $\bar{e}_{j}$.
Choose for each $j$ a lift $c_{j}$ of $\bar{c}_{j}$ (with $c_{0}$
the constant path at $x$).  For fixed $j$ choose $f,g,h \in G$
such that $g\cdot d_{j}$ has the same initial point as $c_{j-1}$
and the sums \[ f\cdot c_{j} + g\cdot d_{j},\quad h\cdot e_{j} +
c_{j-1} \] are defined.  This is possible because of the boundary
relations between the projections in $X/G$ of the various paths.

Now our assumption (ii:b) of Conditions  \ref{P:9.9.7}  implies
that there is an element $k \in G$ such that the following paths
in $X$
\[ k\cdot (f\cdot c_{j} + g\cdot d_{j}), \quad h\cdot e_{j} +
c_{j-1}
\] are homotopic rel end points in $X$.  On applying $\phi$ to
homotopy classes of paths in $X$ and using equations such as $\phi
(g\cdot \gamma) = \phi (\gamma)$ we find that \begin{align*} \phi
[e_{j}] + \phi
[c_{j-1}] &= \phi [h\cdot e_{j}] + \phi [c_{j-1}]\\
        &= \phi [h\cdot e_{j} + c_{j-1}]\\
        &= \phi [k\cdot (f\cdot c_{j} + g\cdot d_{j})] \\% page 9.9-10
        &= \phi [k\cdot f\cdot c_{j}] + \phi [g\cdot d_{j}]\\
        &= \phi [c_{j}] + \phi [d_{j}]. \end{align*} This proves that
\[ \phi [e_{j}] = \phi [c_{j}] + \phi [d_{j}] - \phi [c_{j-1}].
\] It follows easily that  \[ \phi [a_{i-1}] = \phi [a_{i}]. \]
and hence by induction on $i$ that $\phi [a] = \phi [b]$.  From
this it follows that $\phi^{*} : \pi (X/G) \to \Phi$ is a well
defined function such that $\phi^{*}p = \phi$.  The uniqueness of
$\phi^{*}$ is clear since $p_{*}$ is surjective on elements, by
the path lifting property of Conditions \ref{P:9.9.7}.  The proof
that $\phi^{*}$ is a morphism is simple.  This completes the proof
of Proposition \ref{P:9.9.8}. \end{pf}

In the next section, we introduce some further constructions in
the theory of groupoids and groups acting on groupoids, in order
to interpret Proposition \ref{P:9.9.8} in a manner suitable for
calculations. Once again, we will find that an apparently abstract
result involving a universal property can, when appropriately
interpreted, lead to specific calculations.

\begin{center}
EXERCISES 1
\end{center}
\begin{enumerate}[\hspace{-0.95em} 1.]                              \item
\label{Ex:9.9:1} Let $\lambda \in \Real$ and let the additive
group $\Real$ of real numbers act on the torus $T = \ES^{1} \times
\ES^{1}$ by \[ t\cdot (e^{2\pi i\theta},e^{2\pi i\phi}) = (e^{2\pi
i(\theta+t)}, e^{2\pi i(\phi + \lambda t)}) \] for $t,\theta,\phi
\in \Real$. Prove that the orbit space has the indiscrete topology
if and only if $\lambda$ is irrational.  [You may assume that the
group generated by 1 and $\lambda$ is dense in $\Real$ if and only
if $\lambda$ is irrational.]

\item \label{Ex:9.9:2} Let $G$ be a group and let $X$ be a
$G$-space.  Prove that the quotient map $p : X \to X/G$ has the
following universal property: if $Y$ is a space and $f : X \to Y$
is a map such that $f(g \cdot x) = f(x)$ for all $x \in X$ and $g
\in G$, then there is a unique map $f^{*} : X/G \to Y$ such that
$f^{*}p = f$.

\item \label{Ex:9.9:3} Let $X,Y,Z$ be $G$-spaces and let $f : X
\to Z, h : Y \to Z$ be $G$-maps ({\it i.e.}  $f(g \cdot x) =
g\cdot f(x)$ for all $g \in G$ and $x \in X$, and similarly for
$h$).  Let $W = X \times_{Z} Y$ be the pullback.  Prove that $W$
becomes a $G$-space by the action $g\cdot (x,y) = (g \cdot x,g
\cdot y)$.  Prove also that if $Z = X/G$ and $f$ is the quotient
map, then $W/G$ is homeomorphic to $Y$.

\item \label{Ex:9.9:4} Let $G$ and $\Gamma$ be groupoids and let
$w : \Gamma \to \Ob(G)$ be a morphism where $\Ob(G)$ is considered
as a groupoid with identities only.  An {\em action of $G$ on
$\Gamma$ via $w$} is an assignment to each $g \in
G(x,y)$ and % page 9.9-11
$\gamma \in w^{-1}[x]$ an element
$g\cdot \gamma \in w^{-1}[y]$ and with the usual rules: $h\cdot
(g\cdot \gamma) = (hg)\cdot \gamma ; 1\cdot \gamma = \gamma ;
g\cdot (\gamma + \delta) = g\cdot \gamma + g\cdot \delta$.  In
this case $\Gamma$ is called a $G$-groupoid.  Show how to define a
category of $G$-groupoids so that this category is equivalent to
the functor category $\Fun(G,\Set)$.

\item \label{Ex:9.9:5} Prove that if $\Gamma$ is a $G$-groupoid
via $w$, then $\pi_{0}\Gamma$ becomes a $G$-set via $\pi_{0}(w)$.

\item \label{Ex:9.9:6} If $\Gamma$ is a $G$-groupoid via $w$, then
the action is {\em trivial} if for all $x,y \in \Ob(G), g,h \in
G(x,y)$ and $\gamma \in  w^{-1}[x]$, we have $g\cdot \gamma =
h\cdot \gamma$.  Prove that the action is trivial if for all $x
\in  \Ob(G)$, the action of the group $G(x)$ on the groupoid
$w^{-1}[x]$ is trivial.  Prove also that $\Gamma$ contains a
unique maximal subgroupoid $\Gamma^{G}$ on which $G$ acts
trivially.  Give examples to show that $\Gamma^{G}$ may be empty.

\item \label{Ex:9.9:7} Continuing the previous exercise, define a
{\em $G$-section} of $w$ to be a morphism $s : \Ob(G) \to \Gamma$
of groupoids such that $ws = 1$ and $s$ commutes with the action
of $G$, where $G$ acts on $\Ob(G)$ via the source map by $g \cdot
x = y$ for $g \in G(x,y)$.  Prove that $\Gamma^{G}$ is non-empty
if $w$ has a $G$-section, and that the converse holds if $G$ is
connected.  Given a $G$-section $s$, let $\Gamma^{G}(s)$ be the
set of functions $a : \Ob(G) \to \Gamma$ such that $wa = 1$ and a
commutes with the action of $G$ (but we do not assume $a$ is a
morphism).  Show that $\Gamma^{G}(s)$ forms a group under addition
of values, and that if $G$ is connected and $x \in \Ob(G)$, then
$\Gamma^{G}(s)$ is isomorphic to the group $\Gamma (sx)^{G(x)}$ of
fixed points of $\Gamma (sx)$ under the action of $G(x)$.
\end{enumerate}

% page 9.10-1
\section{The computation of orbit groupoids:
quotients and semidirect products}  \label{S:9.10}

The theory of quotient groupoids is modelled on that of quotient
groups, but differs from it in important respects. In particular,
the First Isomorphism Theorem of group theory (that every
surjective morphism of groups is obtained essentially by factoring
out its kernel) is no longer true for groupoids, so we need to
characterise those groupoid morphisms (called quotient morphisms)
for which this isomorphism theorem holds. The next two
propositions achieve this; they were first proved in
\cite{Higgins1}.

Let $f : K \to H$ be a morphism of groupoids.  Then $f$ is said to
be a {\em quotient morphism} if $\Ob(f) : \Ob(K) \to \Ob(H)$ is
surjective and for all $x,y$ in $\Ob(K), f : K(x,y) \to H(fx,fy)$
is also surjective.  Briefly, we say $f$ is {\em object
surjective} and {\em full.}

\begin{Prop} \label{P:9.10.1} Let $f : K \to H$ be a quotient
morphism of groupoids.  Let $N = \Ker f$.  The following hold:
\begin{enumerate}[(i)] \item \label{P:9.10.1:i} If $k,k'  \in  K$,
then $f(k) = f(k')$ if and only if there are elements $m,n \in N$
such that $k'  = m + k + n$.

\item \label{P:9.10.1:ii} If $x$  is an object of $K$, then
$H(fx)$  is isomorphic to the quotient group $K(x)/N(x)$.
\end{enumerate} \end{Prop}

\begin{pf} (\ref{P:9.10.1:i}) If $k,k'$ satisfy $k'  = m + k + n$
where $m,n \in N$, then clearly $f(k) = f(k')$.

Suppose conversely that $f(k) = f(k')$, where $k \in K(x,y), k'
\in K(x', y')$.  Then $$f(x) = f(x'), f(y) = f(y').$$ $$
\spreaddiagramrows{3.0pc}  \spreaddiagramcolumns{3.0pc}
\def\labelstyle{\textstyle}\diagram  x \rto^k & y \ddashed^m
|>\tip \\ x' \udashed^n |>\tip \rto_{k'} & y'  \enddiagram $$
Since $f : K(y,y') \to H(fy,fy')$ is surjective, there is an
element $m \in K(y,y')$ such that $f(m) = 0_{f(y)}$.  Similarly,
there is an element $n \in K(x' ,x)$ such that $f(n) = 0_{f(x)}$.
It follows that if \[ n'  = -k'  + m + k + n  \in K(x'), \] then
$f(n') = 0_{f(x')}$, and so $n' \in N$.  Hence $k' = m+k+n-n'$,
where $m,n - n' \in N$.  This proves (\ref{P:9.10.1:i}).

(\ref{P:9.10.1:ii}) By definition of quotient morphism, the
restriction $f' : K(x) \to H(fx)$ is surjective.  Also by
\ref{P:9.10.1}(i), if $f' (k) = f' (k')$ for $k,k' \in K(x)$, then
there are $m,n \in N(x)$ such that $k' = m + k + n$.  Since $N(x)$
is normal in $K(x)$, there is an $m'$ in $N(x)$ such that $m + k =
k + m'$.  Hence $k' + N(x) = k + N(x)$.  Conversely, if $k' + N(x)
= k + N(x)$ then $f' (k') = f' (k)$. So $f'$ determines an
isomorphism $K(x)/N(x) \to H(fx)$. \end{pf}

We recall the definition of normal subgroupoid.

Let $G$ be a groupoid.  A subgroupoid $N$ of $G$ is called {\em
normal} if $N$ is wide in $G$ (i.e. $\Ob(N)= \Ob(G)$) and, for any
objects $x, y$ of $G$ and $a \in G (x, y)$, $aN(x)a^{-1} \subseteq
N(y)$, from which it easily follows that $$aN(x)a^{-1} = N(y).$$

We now prove a converse of the previous result.  That  is, we
% page 9.10-2
suppose given a normal subgroupoid $N$ of a groupoid $K$ and use
(i) of Proposition \ref{P:9.10.1} as a model for constructing a
quotient morphism $p : K \to K/N$.

The object set of $K/N$ is to be $\pi_{0}N$, the set of components
of $N$.  Recall that a normal subgroupoid is, by definition, wide
in $K$, so that $\pi_{0}N$ is also a quotient set of $X = \Ob(K)$.
Define a relation on the elements of $K$ by $k' \sim {}k$ if and
only if there are elements $m,n$ in $N$ such that $m + k + n$ is
defined and equal to $k'$.  It is easily checked, using the fact
that $N$ is a subgroupoid of $K$, that $\sim {}$ is an equivalence
relation on the elements of $K$. The set of equivalence classes is
written $K/N$.  If $\cls k$ is such an equivalence class, and $k
\in K(x,y)$, then the elements $\cls x$, $\cls y$ in $\pi_{0}N$
are independent of the choice of $k$ in its equivalence class.  So
we can write $\cls k \in K/N(\cls x, \cls y)$. Let $p : K \to K/N$
be the quotient function.  So far, we have not used normality of
$N$.  Not surprisingly, normality is used to give $K/N$ an
addition which makes it into a groupoid.

Suppose  \[ \cls k_{1} \in  (K/N)(\cls x, \cls y),\quad  \cls
k_{2} \in (K/N)(\cls y, \cls z). \]

Then we may assume $k_{1} \in K(x,y), k_{2} \in K(y', z)$, where
$y\~{}y'$ in $\pi_{0}N$.  So there is an element $l \in N(y,y')$,
and we define \[ \cls k_{2} + \cls k_{1} = \cls (k_{2} + l +
k_{1}). \]

We have to show that this addition is well defined.  Suppose then
\begin{align*} k'_{1} &= m_{1} + k_{1} + n_{1},\\ k'_{2} &= m_{2}
+ k_{2} + n_{2}, \end{align*} where $m_{1},n_{1},m_{2},n_{2} \in
N$.  Choose any $l'$ such that $k'_{2} + l'  + k'_{1}$  is
defined.  Then we have the following diagram, in which $a,a'$ are
to be defined:
 $$
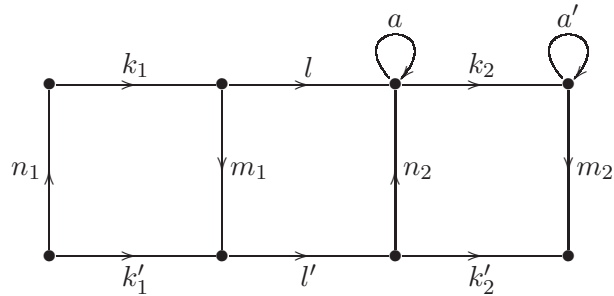
\spreaddiagramrows{3.0pc} \spreaddiagramcolumns{3.0pc}
\objectmargin{0pc} \objectwidth{0pc}
\def\labelstyle{\textstyle}\xymatrix{ \bullet \ar @{-} [r]^{k_1}|\tip
\ar @{-} [d] _{n_1}|{\rotate\tip}&\bullet \ar @{-} [r] ^{l}|\tip
\ar @{-} [d] ^{m_1}|\tip  &\bullet  \ar @{-} [r] ^{k_2}|\tip
\ar@(ul,ur)[]^a
\ar @{-} [d] ^{n_2}|{\rotate\tip}& \ar @{-} [d] ^{m_2}|\tip \bullet  \ar@(ul,ur)[]^{a'}  \\
\bullet \ar @{-} [r] _{k'_1}|\tip&\bullet \ar @{-} [r] _{l'}|\tip
& \bullet \ar @{-} [r] _{k'_2}|\tip &\bullet} $$
  Let $a = n_{2} + l' +
m_{1} - l$.  Then $a \in N$, and $l = -a + n_{2} + l' + m_{1}$.
Since $N$ is normal there is an element $a' \in N$ such that $a' +
k_{2} = k_{2} - a$.  Hence \begin{align*} k_{2} + l +
k_{1} &= k_{2} - a + n_{2} + l' + m_{1} + k{}_{1}\\
      &= a'  + k_{2} + n_{2} + l'  + m_{1} + k_{1} \\% page 9.10-3
      &= a'  - m_{2} + k'_{2} + l'  + k'_{1} - n_{1}.  \end{align*}
Since $a',m_2,n_1 \in N$, we obtain $\cls (k_{2} + l + k_{1}) =
\cls (k'_{2} + l' + k'_{1})$ as was required.

Now we know that the addition on $K/N$ is well defined, it is easy
to prove that the addition is associative, has identities, and has
inverses.  We leave the details to the reader. So we know that
$K/N$ becomes a groupoid.

\begin{Prop} \label{P:9.10.2} Let $N$ be a normal subgroupoid of
the groupoid $K$, and let $K/N$ be the groupoid just defined. Then
\begin{enumerate}[(i)] \item \label{P:9.10.2:i}  the quotient function
$p : k \mapsto \cls k$ is a quotient morphism $K \to K/N$ of
groupoids;

\item \label{P:9.10.2:ii} if $f : K \to H$  is any morphism of
groupoids such that Ker$ f$ contains $N$, then there is a unique
morphism $f^{*} : K/N \to H$ such that $f^{*}p = f$.
\end{enumerate} \end{Prop}

\begin{pf} The proof of \eqref{P:9.10.2:i} is clear.  Suppose $f$
is given as in \eqref{P:9.10.2:ii}.  If $m + k + n$ is defined in
$K$ and $m,n \in  N$, then $f(m + k + n) = f(k)$.  Hence $f^{*}$
is well defined on $K/N$ by $f^{*}(\cls k) = f(k)$.  Clearly
$f^{*}p = f$.  Since $p$ is surjective on objects and elements,
$f^{*}$ is the only such morphism. \end{pf}

In order to apply these results, we need generalisations of some
facts on normal closures which were given in Section 3 of Chapter
8 for the case of a family $R(x)$ of subsets of the object groups
$K(x), x \in \Ob(K)$, of a groupoid $K$.  The argument here is
based on \cite[Exercise 4, p.95]{Higgins4}.

Suppose  that $R$ is {\em any} set of elements of the groupoid
$K$. The {\em normal closure} of $R$ in $K$ is the smallest normal
subgroupoid $N(R)$ of $K$ containing $R$. Clearly $N(R)$ is the
intersection of all normal subgroupoids of $K$ containing $R$, but
it is also convenient to have an explicit description of $N(R)$.

\begin{Prop} \label{P:9.10.normal}
Let $\lan R \ran$ be the wide subgroupoid of $K$ generated by $R$.
Then the normal closure $N(R)$ of $R$ is the subgroupoid of $K$
generated by $\lan R \ran$ and all conjugates $k h k \io $ for $k
\in K, h \in \lan R \ran$.
\end{Prop}
\begin{pf}
Let $\widehat{R}$ be the subgroupoid of $K$ generated by $\lan R
\ran$ and all conjugates $k h k \io $ for $k \in K, h \in \lan R
\ran$. Clearly any normal subgroupoid of $K$ containing $R$
contains $\widehat{R}$, so it is sufficient to prove that
$\widehat{R}$ is normal.

Suppose then that $k+a-k$ is defined where $k \in K$ and  $a \in
\widehat{R}$ so that
\[ a = r_{1} + c_{1} + r_{2} + c_{2} + \cdots  + r_{l} +
c_{l} + r_{l+1} \] where each $r_i \in \lan R \ran$ and each
$c_i=k_i +h_i -k_i$ is a conjugate of a loop $h_i$ in $\lan R
\ran$ by an element $k_i \in K$.
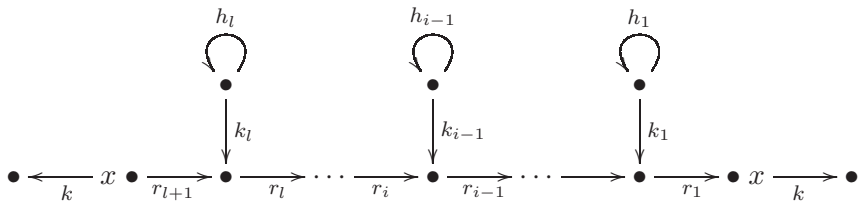
$$\xymatrix{&& \bullet \ar [d]^{k_l} \ar @(ur,ul) _{h_l}
&& \bullet \ar [d]^{k_{i-1}} \ar @(ur,ul) _{h_{i-1}} &&
\bullet \ar [d]^{k_1} \ar @(ur,ul) _{h_1}& &\\
\bullet  & \ar [l]^k x \;\bullet \ar[r]_{r_{l+1}}& \bullet
\ar[r]_{r_l}& \cdots   \ar [r]_{r_i}   & \bullet \ar [r]_{r_{i-1}}
& \cdots \ar[r] & \bullet \ar [r] _ {r_1}& \bullet \; x \ar [r] _k
& \bullet }$$ Then $a$ is a loop, since $k+a-k$ is defined, and so
also is
$$b=r_1+ r_2+ \cdots +r_{l+1}.$$ Let
$$d_i=k+r_1+\cdots+r_i+c_i-r_i -\cdots -r_1-k$$ so that $d_i$
is a conjugate of a loop in $\lan R \ran$ for $i=1, \ldots, l$.
Then it is easily checked that
$$k+a-k= d_1+\cdots + d_l+k+b-k$$
and hence $k+a-k \in\widehat{R}$.
\end{pf}
Notice that the loop $b$ in the proof belongs to $\lan R \ran $
rather than to $R$, and this shows why it is not enough just to
take $N(R)$ to be the subgroupoid generated by $R$ and conjugates
of loops in $R$.

The elements of $N(R)$ as constructed above may be called the {\it
consequences} of $R$.

We next give the definition of the {\it semidirect product} of a
group with a groupoid on which it acts.  Let $G$ be a group and
let $\Gamma$ be a groupoid with $G$ acting on the left.  The {\it
semidirect product groupoid} $\Gamma \rtimes G$ has object set
$\Ob(\Gamma)$ and arrows $x \to y$ the set of pairs $(\gamma, g)$
such that $g \in G$ and $\gamma \in \Gamma (g \cdot x,y)$.  The
sum of $(\gamma ,g) : x \to y$ and $(\delta, h) : y \to z$ in
$\Gamma \rtimes G$ is defined to be \[ (\delta, h) + (\gamma, g) =
(\delta + h\cdot \gamma, hg). \] This is easily remembered from
the following picture. $$
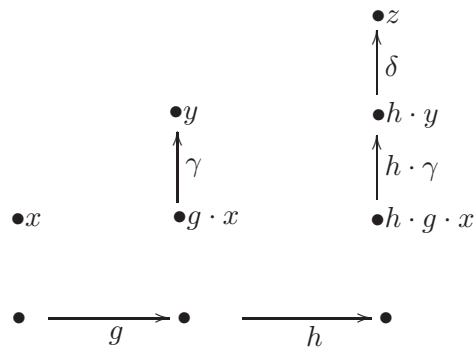
\def\labelstyle{\textstyle}  \spreaddiagramcolumns{1.5pc}
%\spreaddiagramrows{2pc}
\diagram &&\bullet z \hspace{2.55em}\\
& \bullet y \hspace{1.55em} & \bullet h\cdot y
\hspace{1.2em}\uto<3.5ex>_{\delta} \\ \bullet x & \bullet g \cdot
x \uto<2.2ex>_{\gamma} & \bullet h \cdot g\cdot x\uto<3.5ex>_{h
\cdot\gamma} \\ \bullet \hspace{0.5em}\rto_g & \bullet
\hspace{1.5em}\rto_h & \bullet \hspace{2.5em} \enddiagram $$
%this figure might be better done in picture environment
%\begin{center} {\bf Fig. 8} \end{center}
\begin{Prop}
\label{P:9.10.3} The above addition makes $\Gamma \rtimes G$ into
a groupoid and the projection $$q : \Gamma \rtimes G \to G, \quad
(\gamma, g) \mapsto g, $$ is a fibration of groupoids. Further:
\begin{enumerate}[(i)]  \item \label{P:9.10.3:i} $q$  is a quotient
morphism if and only if $\Gamma$  is connected; \item
\label{P:9.10.3:ii} $q$  is a covering morphism if and only if
$\Gamma$  is discrete;  \item \label{P:9.10.3:iii} $q$ maps
$(\Gamma \rtimes G)(x)$ isomorphically to $G$ for all $x \in
\Ob(\Gamma)$ if and only if $\Gamma$ has trivial object groups and
$G$ acts trivially on $\pi_{0}\Gamma$. \end{enumerate}
\end{Prop}

\begin{pf} The proof of the axioms for a groupoid is easy, the
negative of $(\gamma, g)$ being $(g^{-1}\cdot (-\gamma),g^{-1})$.
We leave the reader to check associativity.

To prove that $q$ is a fibration, let $g \in G$ and $x \in
\Ob(\Gamma)$.  Then $(0_{g \cdot x}, g)$ has source $x$ and maps
by $q$ to $g$.

We now prove (\ref{P:9.10.3:i}).  Let $x,y$ be objects of
$\Gamma$. Suppose $q$ is a quotient morphism.  Then $q$ maps
$(\Gamma \rtimes G)(x,y)$ surjectively to $G$ and so there is an
element $(\gamma, g)$ such that $q(\gamma, g) = 1$.  So $g = 1$
and $\gamma \in \Gamma (x,y)$.  This proves $\Gamma$ is connected.

Suppose $\Gamma$ is connected.  Let $g \in G$.  Then there is a
$\gamma \in \Gamma (g \cdot x,y)$, and so $q(\gamma, g) = g$.
Hence $q$ is a quotient morphism.

We now prove \eqref{P:9.10.3:ii}.  Suppose $\Gamma$ is discrete,
so that $\Gamma$ may be thought of as a set on which $G$ acts.
Then $\Gamma \rtimes G$ is simply the % page 9.10-5
covering groupoid of the action as constructed in a previous
section. So $q$ is a covering morphism.

Let $x$ be an object of $\Gamma$.  If $\gamma$ is an element of
$\Gamma$ with source $x$ then $(\gamma, 1)$ is an element of
$\Gamma \rtimes G$ with source $x$ and which lifts $1$.  So if $q$
is a covering morphism then the star of $\Gamma$ at any $x$ is a
singleton, and so $\Gamma$ is discrete.

The proof of \eqref{P:9.10.3:iii} is best handled by considering
the exact sequence based at $x \in \Ob(\Gamma)$ of the fibration
$q$. This exact sequence is by 7.2.10 of \cite{Brown:1988}\[ 1 \to
\Gamma (x) \to (\Gamma \rtimes G)(x)\labto{q'} G \to \pi_{0}\Gamma
\to \pi_{0}(\Gamma \rtimes G) \to 1. \] It follows that $q'$ is
injective if and only if $\Gamma (x)$ is trivial. Exactness also
shows that $q'$ is surjective for all $x$ if and only if the
action of $G$ on $\pi_{0} \Gamma$ is trivial. \end{pf}

Here is a simple application of the definition of semidirect
product which will be used later.

\begin{Prop} \label{P:9.10.4} Let $G$ be a group and let
$\Gamma$ be a $G$-groupoid.  Then the formula $$(\gamma, g)\cdot
\delta = \gamma + g\cdot \delta $$  for $\gamma, \delta \in
\Gamma, g \in G,$ defines an action of $\Gamma \rtimes G$ on the
set $\Gamma$ via the target map $\tau : \Gamma \to \Ob(\Gamma)$.
\end{Prop}

\begin{pf} This says in the first place that if $(\gamma, g) \in
(\Gamma \rtimes G)(y,z)$ and $\delta$ has target $y$, then $\gamma
+ g\cdot \delta$ has target $z$, as is easily verified. The axioms
for an action are easily verified.  The formula for the action
also makes sense if one notes that \[ (\gamma, g)(\delta, 1) =
(\gamma  + g\cdot \delta, g). \] \end{pf}

If $X$ is a $G$-space, and $x \in X$, let $\sigma (X,x,G)$ be the
object group of the semidirect product groupoid $\pi X \rtimes G$
at the object $x$.  This group is called by Rhodes in
\cite{Rhodes1} and \cite{Rhodes2} the {\em fundamental group of
the transformation group} (although he defines it directly in
terms of paths).  The following result from \cite{Rhodes2} gives
one of the reasons for its introduction.

\begin{corollary} \label{Cor:9.10.4:1} If $X$ is a $G$-space, $x
\in X$, and the universal cover $\tilde{X}_{x}$ exists, then the
group $\sigma (X,x,G)$ has a canonical action on $\tilde{X}_{x}$.
\end{corollary}

\begin{pf} By 9.5.8 of \cite{Brown:1988}, we may identify the
universal cover $\tilde{X}_{x}$ of $X$ at $x$ with $\St_{\pi X}x$.
The function $\pi X \to \pi X, \delta \mapsto -\delta$, transports
the action of $\pi X \rtimes G$ on $\pi X$ via the target map
$\tau$ to an action of the same groupoid on $\pi X$ via the source
map $\sigma$.  Hence the object group $(\pi X \rtimes G)(x)$ acts
on $\St_{\pi X}x$ by \[ (\gamma,g)*\delta = -((\gamma,g)\cdot
(-\delta))=g\cdot \gamma-\gamma. \] The continuity of the action
follows easily
from the detailed % page 9.10-6
description of the lifted topology (see also the remarks on
topological groupoids after 9.5.8 of \cite{Brown:1988}). \end{pf}
Now we start using the semidirect product to compute orbit
groupoids. The next two results may be found in
\cite{HigginsTaylor1}, \cite{JTaylor1} and \cite{JTaylor2}.

\begin{Prop} \label{P:9.10.5} Let $N$ be the normal closure in
$\Gamma \rtimes G$ of the set of elements of the form $(0_{x},g)$
for all $x \in \Ob(\Gamma)$ and $g \in G$ .  Let $p$ be  the
composite  \[ \Gamma \labto{i} \Gamma \rtimes G \labto{\nu}
(\Gamma \rtimes G)/N, \] in which the first morphism is $\gamma
\mapsto (\gamma, 1)$ and the second morphism is the quotient
morphism. Then
\begin{enumerate}[(i)] \item \label{P:9.10.5:i} $p$  is a surjective
fibration; \item \label{P:9.10.5:ii} $p$  is an orbit morphism and
so determines an isomorphism \[ \Gamma \orbitgpd G \cong  (\Gamma
\rtimes G)/N; \] \item \label{P:9.10.5:iii} the function
$\Ob(\Gamma) \to \Ob(\Gamma \orbitgpd G)$  is an orbit map, so
that $\Ob(\Gamma \orbitgpd G)$  may be identified with the orbit
set $\Ob(\Gamma)/G$. \end{enumerate} \end{Prop}

\begin{pf} Let $\Delta = (\Gamma \rtimes G)/N$.  We first derive
some simple consequences of the definition of $\Delta$.  Let
$\gamma \in \Gamma(x,y), g,h \in  G$.  Then \begin{equation}
\label{Eq:9.10:1} (0_{g \cdot y},g) + (\gamma, 1) = (g\cdot
\gamma, g), \end{equation} \begin{equation} \label{Eq:9.10:2}
(\gamma, 1) + (0_{x},h) = (\gamma, h). \end{equation} It follows
that in $\Delta$ we have  \begin{equation} \label{Eq:9.10:3} \nu
(h\cdot \gamma, 1) = \nu (\gamma, g). \end{equation} Note also
that the set $R$ of elements of $\Gamma \rtimes G$ of the form
$(0_{g \cdot x},g)$ is a subgroupoid of $G$, since  \[ (0_{hg
\cdot x},h) + (0_{g \cdot x},g) = (0_{hg \cdot x},hg), \] and
$-(0_{g \cdot x},g) = (0_{x},g^{-1})$.  It follows that $\pi_{0}N
= \pi_{0}R = \Ob(\Gamma)/G$, the set of orbits of the action of
$G$ on $\Ob(\Gamma)$.  Hence $\Ob(p)$ is surjective. This proves
\eqref{P:9.10.5:iii}, once we have proved \eqref{P:9.10.5:ii}.

We now prove easily that $p : \Gamma \to \Delta$ is a fibration.
 Let $(\gamma, g) : x \to y$ in $\Gamma \rtimes G$ be a
representative of an element of $\Delta$, and suppose $\nu z = \nu
x$, where $z \in \Ob(\Gamma)$.  Then $z$ and $x$ belong to the
same orbit and so there is an element $h$ in $G$ such that $h
\cdot x = z$.  Clearly $h\cdot \gamma$ has source $z$ and by
\eqref{Eq:9.10:3}, $p(h\cdot \gamma) = \nu(\gamma,g)$.

Suppose now $g \in G$ and $\gamma : x \to y$ in $\Gamma$.  Then by
\eqref{Eq:9.10:3} $p(g\cdot \gamma) = p(\gamma)$.  This verifies
\eqref{P:9.9.6:i}.

To prove the other condition for an orbit morphism, namely
Definition \ref{P:9.9.6}(ii), suppose $\phi : \Gamma \to \Phi$ is
a morphism of groupoids such that $\Phi$ has a trivial action of
the group $G$ and $\phi (g\cdot \gamma ) = \phi (\gamma)$ for all
$\gamma \in \Gamma$ and $g \in G$.  Define $\phi' : \Gamma \rtimes
G \to \Phi$ on objects by $\Ob(\phi)$ and on elements by $(\gamma,
g) \mapsto \phi (\gamma)$.  That $\phi$ is a morphism follows from
the trivial action of $G$ on $\Phi$, since
\begin{align*} \phi' ((\delta, h) + (\gamma, g)) &= \phi (\gamma
 + h.\delta)\\
&= \phi (\gamma) + \phi (h.\delta)\\
&= \phi (\gamma) + \phi(\delta)\\
&= \phi' (\delta, h) + \phi' (\gamma, g). \end{align*} Also $\phi'
(0_{x},g) = \phi (0_{x}) = 0_{\phi x}$, and so $N \subseteq \Ker
\phi'$.  By \eqref{P:9.10.2:ii}, there is a unique morphism
$\phi^{*} : (\Gamma \rtimes G)/N \to \Phi$ such that $\phi^{*}\nu
= \phi'$.  It follows that $\phi^{*}p = \phi^{*}\nu i = \phi'i =
\phi$.  The uniqueness of $\phi^{*}$ follows from the fact that
$p$ is surjective on objects and on elements.

Finally, the isomorphism $\Gamma \orbitgpd G \cong \Delta$ follows
from the universal property. \end{pf}  In order to use the last
result  we analyse the morphism $p : \Gamma \to \Gamma \orbitgpd
G$ in some special cases. The construction of the orbit groupoid T
in Proposition \ref{P:9.10.5} is what makes this possible.

\begin{Prop} \label{P:9.10.6} The orbit morphism $p : \Gamma \to
\Gamma \orbitgpd G$ is a fibration whose kernel is generated as a
subgroupoid of $\Gamma$ by all elements of the form $\gamma -
g\cdot \gamma$ where $g$ stabilises the source of $\gamma$.
Furthermore, \begin{enumerate}[(i)] \item \label{P:9.10.6:i} if
$G$ acts freely on $\Gamma$, by which we mean no non-identity
element of $G$ fixes an object of $\Gamma$, then $p$ is a covering
morphism;

\item \label{P:9.10.6:ii} if $\Gamma$  is connected and $G$  is
generated by those of its elements which fix some object of
$\Gamma$, then $p$ is a quotient morphism; in particular, $p$ is a
quotient morphism if the action of $G$ on $\Ob(\Gamma)$ has a
fixed point;

\item \label{P:9.10.6:iii} if $\Gamma$  is a tree groupoid, then
each object group of $\Gamma \orbitgpd G$ is isomorphic to the
factor group of $G$ by the (normal) subgroup of $G$ generated by
elements which have fixed points. \end{enumerate} \end{Prop}

\begin{pf} We use the description of $p$ given in the previous
Proposition \ref{P:9.10.5}, which already implies that $p$ is a
fibration.

Let $R$ be the subgroupoid of $\Gamma \rtimes G$ consisting of
elements $(0_{g \cdot x},g), g \in G$.  Let $N$ be the normal
closure of $R$.  By the construction of the normal closure in
Proposition \ref{P:9.10.normal}, the elements of $N$ are sums of
elements of $ R $ and conjugates of loops in $R$ by elements of
$\Gamma \rtimes G$. So let $(0_{g \cdot x},g)$ be a loop in $R$.
Then $g \cdot x = x$. Let $(\gamma, h) : x \to y$ in $\Gamma
\rtimes G$, so that $\gamma : h \cdot x \to y$.  Then we check
that \[ (\gamma, h) + (0_{x},g) - (\gamma, h)
= (\gamma  - hgh^{-1}\cdot\gamma, hgh^{-1}). \] % page 9.10-8
Writing $k = hgh^{-1}$, we see that $(\gamma - k\cdot \gamma, k)
\in  N$ if $k$ stabilises the initial point of $\gamma$.

Now $\gamma \in \Ker p$ if and only if $(\gamma, 1) \in N$.
Further, if $(\gamma, 1) \in N$ then $(\gamma, 1)$ is a
consequence of $R$ and so $(\gamma, 1)$ is equal to \[ (\gamma_{1}
-k_{1}\cdot \gamma_{1},k_{1}) + (\gamma_{2} -k_{2}\cdot
\gamma_{2},k_{2}) + \dots + (\gamma_{r} -k_{r}\cdot
\gamma_{r},k_{r}) \] for some $\gamma_{i}, k_{i}$ where $k_{i}$
stabilises the initial point of $\gamma_{i}, i = 1,\dots, r$. Let
$h_{1} = 1, h_{i} = k_{1}\dots k_{i-1}(i\geqslant 2), \delta_{i} =
h_{i}.\gamma, g_{i} = h_{i}k_{i}h^{-1}_{i}, i\geqslant 1$.  Then
\[ (\gamma, 1) = (\delta_{1} -g_{1}\cdot \delta_{1} +
\delta_{2}g_{2}\cdot \delta_{2} + \dots + \delta_{r} -g_{r}\cdot
\delta _{r},1) \] and so $\gamma$ is a sum of elements of the form
$\delta - g\cdot \delta$ where $g$ stabilises the initial point of
$\delta$.  This proves our first assertion.

The proof of \eqref{P:9.10.6:i} is simple.  We know already that
$p : \Gamma \to \Gamma \orbitgpd G$ is a fibration.  If $G$ acts
freely, then by the result just proved, $p$ has discrete kernel.
 It follows that if $x \in \Ob(\Gamma)$, then $p : \St_{\Gamma}x
\to \St_{\Gamma \orbitgpd G}px$ is injective.  Hence $p$ is a
covering morphism.

Now suppose $\Gamma$ is connected and $G$ is generated by those of
its elements which fix some object of $\Gamma$.  To prove $p$ a
quotient morphism we have to show that for $x,y \in \Ob(\Gamma)$,
the restriction $p' : \Gamma (x,y) \to (\Gamma \orbitgpd
G)(px,py)$ is surjective.

Let $(\gamma, g)$ be an element $ x \to y$ in $\Gamma \rtimes G$,
so that $\gamma : g \cdot x \to y$ in $\Gamma$.  Using the
notation of \ref{P:9.10.5}, we have to find $\delta \in \Gamma
(x,y)$ such that $p\delta = \nu (\gamma, g)$.  As shown in
\ref{P:9.10.5}, $\nu (\gamma, g) = \nu (\gamma, 1)$.  By
assumption, $g = g_{n}g_{n-1}\dots g_{1}$ where $g_{i}$ stabilises
an object $x_{i}$, say.  Since $\Gamma$ is connected, there are
elements
\[ \delta_{1} \in  \Gamma (x_{1},x), \;\; \delta_{i} \in  \Gamma
(x_{i}, (g_{i-1}\dots g_{1}) \cdot x), i\geqslant 1. \] The
situation
is illustrated below for $n = 2$. %\ref{Fig:9.8}.
$$
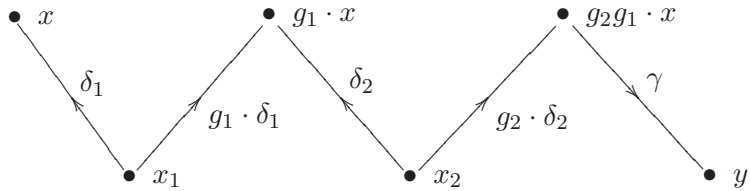
\spreaddiagramrows{2.5pc}  %\spreaddiagramcolumns{2.5pc}
\objectmargin{0pc} \objectwidth{0pc}
\def\labelstyle{\textstyle}
\xymatrix{ **[r]   \bullet  \; \; x&& **[r] \bullet\; \; g_1 \cdot
x
&& **[r] \bullet \;\;  g_2g_1 \cdot x \ar @{-} [dr]| \tip  ^-{\gamma} & \\
& **[r] \bullet \; \; x_1 \ar @{-} [ul]|\tip _{\delta_1}  \ar @{-}
[ur]|\tip _{g_1 \cdot\delta_1} &&**[r] \bullet\; \; x_2\ar @{-}
[ul]|\tip _{\delta_2} \ar @{-} [ur]|\tip _{g_2 \cdot\delta_2} &&
**[r] \bullet \;\;  y}$$

%\begin{center} Figure 4 \end{center}
%\begin{figure} \caption{} \label{Fig:9.8} \end{figure}
Let
$$\delta = \gamma + (g_{n}\cdot \delta_{n} - \delta_{n}) + \cdots  +
(g_{1}\cdot \delta_{1} - \delta_{1}): x \to y.$$  Then $p(\delta)
= \nu (\gamma, 1)$.  This proves \eqref{P:9.10.6:ii}.

For the proof of \eqref{P:9.10.6:iii}, let $x$ be an object of
$\Gamma$. Since $\Gamma$ is a tree groupoid, the projection
$(\Gamma \rtimes G)(x) \to G$ is an isomorphism which sends the
element $(0_{x},g)$, where $g \cdot x = x$, to the element $g$.
Also if $g$ fixes $x$ and $h \in G$ then $hgh^{-1}$ fixes $h \cdot
x$. Thus the image of $N(x)$ is the subgroup $K$ of $G$
generated by elements of $G$ with a % page 9.10-9
fixed point, and $K$ is normal in $G$.  Let $\bar{x}$ denote the
orbit of $x$. By \eqref{P:9.10.1:ii}, and Proposition
\ref{P:9.10.5}, the group $(\Gamma \orbitgpd G)(\bar{x})$ is
isomorphic to the quotient of $(\Gamma \rtimes G)(x)$ by $N(x)$.
Hence $(\Gamma \orbitgpd G)(\bar{x})$ is also isomorphic to $G/K$.
\end{pf}

In Proposition \ref{P:9.10.6}, the result (\ref{P:9.10.6:i})
relates the work on orbit groupoids to work on covering morphisms.
The result (\ref{P:9.10.6:ii}) will be used below.  The result
(\ref{P:9.10.6:iii}) is particularly useful in work on
discontinuous actions on Euclidean or hyperbolic space.  In the
case of a discontinuous action of a group $G$ on a space $X$ which
satisfies the additional condition (\ref{P:9.9.7:ii}) of
\ref{P:9.9.7}, and which has a fixed point $x$, we obtain from
Proposition \ref{P:9.10.6} a convenient description of $\pi
(X/G,\bar{x})$ as a quotient of $\pi (X,x)$. In the more general
case, $\pi (X/G,\bar{x})$ has to be computed as a quotient of
$\sigma (X,x,G) = (\pi X \rtimes G)(x)$, as given in Proposition
\ref{P:9.10.5}.  Actually the case corresponding to
(\ref{P:9.10.6:iii}) was the first to be discovered, for the case
of simplicial actions \cite{Armstrong1}.  The general case
followed from the fact that if $X$ has a universal cover
$\tilde{X}_{x} $ at $X$, then $\sigma (X,x,G)$ acts on
$\tilde{X}_{x}$ with orbit space homeomorphic to $X/G$.

We give next two computations for actions with fixed points.

% These next two propositions are wierd.  Are they supposed to
%be
% examples? \noindent
\begin{example}  \label{P:9.10.7}
Let the group $\Z_{2}= \{1,g\}$ act on the circle $\ES^{1}$ in
which $g$ acts by reflection in the $x$-axis.  The orbit space of
the action can be identified with $E^{1}_{+}$, which is
contractible, and so has trivial fundamental group.

To see how this agrees with the previous results, let $\Gamma =
\pi \ES^{1}, G = \Z_{2}$.  It follows from Proposition
\ref{P:9.9.8} that the induced morphism $\pi \ES^{1} \to \pi
(\ES^{1}/G)$ is an orbit morphism, and so $\pi (\ES^{1}/G) \cong
\Gamma \orbitgpd G$. By Proposition \ref{P:9.10.5}, $p : \Gamma
\to \Gamma \orbitgpd G$ is a quotient morphism, since the action
has a fixed point $1$ (and also $-1$).
 Let the two elements  $a_\pm  \in \Gamma (1,-1)$ be represented by the paths $[0,1]
\to \ES^{1}, t \mapsto  e^{\pm i\pi \,t}$ respectively.  The non
trivial element  $g$ of $G$ satisfies $g\cdot a_+ = a_-$. Hence
the kernel of the quotient morphism $p : \Gamma (1) \to (\Gamma
\orbitgpd G)(p1)$ contains the element $-a_- + a_+$.  But this
element generates $\Gamma (1) \cong \Z$.  So we confirm the fact
that $(\Gamma \orbitgpd G)(p1)$ is the trivial group.
\rule{1em}{0ex} \hfill $\Box$\end{example}

% HUH?  Is this ALL theorem??? [Changed]

Before our next result we state and prove a simple group theoretic
result. First let $H$ be a group. It is convenient to write the
group structure on $H$ as multiplication.  The {\em
abelianisation} $H^{\ab}$ of $H$ is formed from $H$ by imposing
the relations $hk = kh$ for all $h,k \in H$.  Equivalently, it is
the quotient of $H$ by the (normal) subgroup generated by all
commutators $hkh^{-1}k^{-1}$, for all $h,k \in H$.

\begin{Prop} The quotient $(H \times H)/K$ of $H \times H$ by
the normal subgroup $K$ of $H \times H$ generated by the elements
$(h,h^{-1})$ is isomorphic to $H^{ab}$. \end{Prop}
\begin{pf} We can regard $H \times H$ as the group with
generators $[h], \<k\>$ for all $h,k \in H$ and relations $[hk] =
[h][k], \<hk\> = \<h\>\<k\>, [h]\<k\> = \<k\>[h]$ for all $h,k \in
H$, where we may identify $[h] = (h,1), \<k\> = (1,k), [h]\<k\> =
(h,k)$. Factoring out by $K$ imposes the additional relations
$[h]\<h^{-1}\> = 1$, or equivalently $[h] = \<h\>$, for all $h \in
H$. It follows that $(H \times H)/K$ is obtained from $H$ by
imposing the additional relations $hk = kh$ for all $h,k \in H$.
 \end{pf}
 \begin{Def}[The symmetric square of a
space] \label{P:9.10.8}   Let $G = \Z_{2}$ be the cyclic group of
order $2$, with non-trivial element $g$. For a space $X$, let $G$
act on the product space $X \times X$ by interchanging the
factors, so that $g\cdot (x,y) = (y,x)$.  The fixed point set of
the action is the diagonal of $X \times X$.  The orbit space is
called the % page 9.10-10
{\em symmetric square} of $X$, and is
written $Q^{2}X$.
\end{Def}
\begin{Prop}
Let $X$ be a connected, Hausdorff, semilocally $1$-connected
space, and let $x \in X$. Let $\< x\>$ denote the class in
$Q^{2}X$ of $(x,x)$. Then the fundamental group $\pi(Q^{2}X,\<
x\>)$ is isomorphic to $\pi(X,x)^\ab$, the fundamental group of
$X$ at $x$ made abelian.
\end{Prop}
\begin{pf}
 Since $G=\Z_2$ is finite, the action is
discontinuous.  Because of the assumptions on $X$, we can apply
Proposition \ref{P:9.9.8}, and hence also the results of this
section. We deduce that $p_{*} : \pi (X) \times \pi( X) \to \pi(
Q^{2}X)$ is a quotient morphism and that if $x \in X, z = p(x,x)$,
then the kernel of the quotient morphism \[ p' : \pi (X,x) \times
\pi (X,x) \to \pi(Q^{2}X,z) \] is the normal subgroup $K$
generated by elements $(a,b) - g\cdot (a,b) = (a-b,b-a), a,b \in
\pi X(y,x)$, for some $y \in X$.  Equivalently, $K$ is the normal
closure of the elements $(c,-c), c \in \pi (X,x)$. The result
follows.
\end{pf}

Taylor \cite{JTaylor2} has given extensions of previous results
which we give without proof.

\begin{Prop} \label{P:9.10.9} Let $G$ be a group and let
$\Gamma$ be a $G$-groupoid.  Let $A$ be a subset of $\Ob(\Gamma)$
such that $A$ is $G$-invariant and for each $g \in G$, $A$ meets
each component of the subgroupoid of $G$ left fixed by the action
of $G$.  Let $\Xi$ be the full subgroupoid of $\Gamma$ on $A$.
Then the orbit groupoid $\Xi \orbitgpd G$ is embedded in $\Gamma
\orbitgpd G$ as the full subgroupoid on the set of objects $A/G$
and the orbit morphism $\Xi \to \Xi \orbitgpd G$ is the
restriction of the orbit morphism $\Gamma \to \Gamma \orbitgpd G$.
\end{Prop}

\begin{corollary} \label{Cor:9.10.9:1} If the action of $G$  on
$X$ satisfies Conditions % page 9.10-11
\ref{P:9.9.7}, and $A$ is a $G$-stable subset of $X$ meeting each
path component of the fixed point set of each element of $G$, then
$\pi (X/G)(A/G)$ is canonically isomorphic to $(\pi XA)\orbitgpd
G$.
\end{corollary}

The results of this section answer some cases of the following
question:

\begin{question} Suppose $p : K \to H$ is a connected covering
morphism, and $x \in \Ob(K)$.  Then $p$ maps $K(x)$ isomorphically
to a subgroup of the group $H(px)$.  What information in addition
to the value of $K(x)$ is needed to reconstruct the group $H(px)$?
\end{question}

There is an exact sequence \[ 0 \to K(x) \to H(px) \to
H(px)/p[K(x)] \to 0 \] in which $H(px)$ and $K(x)$ are groups
while $H(px)/p[K(x)]$ is a pointed set with base point the coset
$p[K(x)]$.  Suppose now that $p$ is a regular covering morphism.
 Then the group $G$ of covering transformations of $p$ is
anti-isomorphic to $H(px)/p[K(x)]$, by \cite[Cor.
9.6.4]{Brown:1988}. Also $G$ acts freely on $K$.

\begin{Prop} \label{P:9.10.10} If $p : K \to H$ is a regular
covering morphism, then $p$ is an orbit morphism with respect to
the action on $K$ of the group $G$ of covering transformations of
$p$.  Hence if $x \in \Ob(K)$, then $H(px)$ is isomorphic to the
object group $(K \rtimes G)(x)$. \end{Prop}

\begin{pf} Let $G$ be the group of covering transformations of
$p$.  Then $G$ acts on $K$.  Let $q : K \to K\orbitgpd G$ be the
orbit morphism.  If $g \in G$ then $pg = p$, and so $G$ may be
considered as acting trivially on $H$.  By Definition
\ref{P:9.9.6}(i), there is a unique morphism $\phi : K\orbitgpd G
\to H$ of groupoids such that $\phi q = p$.  By Proposition
\ref{P:9.10.6}(i), $q$ is a covering morphism.  Hence
 $\phi$ is a covering morphism \cite[9.2.3]{Brown:1988}.  But
$\phi$ is bijective on objects, because $p$ is regular.  Hence
$\phi$ is an isomorphism.

Since $G$ acts freely on $K$, the group $N(x)$ of Proposition
\ref{P:9.10.5} is trivial.  So the description of $H(px)$ follows
from Proposition \ref{P:9.10.5}(ii). \end{pf}

The interest of the above results extends beyond the case where
$G$ is finite, since general discontinuous actions occur in
important applications in complex function theory, concerned with
Fuchsian groups and Kleinian groups.We  refer the reader to
\cite{Beardon1} and \cite{ECD}.

\begin{center}
EXERCISES
\end{center}

\begin{enumerate}[\hspace{-1em} 1.]

\item \label{Ex:9.10:1} Let $f : G \to H$ be a groupoid morphism
with kernel $N$.  Prove that the following are equivalent:
\begin{enumeratein} \item $f$ is a quotient morphism; \item $f$
is surjective and any two vertices of $G$ having the same image in
$H$ lie in the same component of $G$. \end{enumeratein}

\item \label{Ex:9.10:2} Prove that a composite of quotient
morphisms is a quotient morphism.

\item \label{Ex:9.10:3} Let $H$ be a subgroupoid of the groupoid
$G$ with inclusion morphism $i : H \to G$.  Let $f : G \to H$ be a
morphism with kernel $N$. Prove that the following are equivalent:
\begin{enumeratein} \item $f$ is a deformation retraction; \item
$f$ is piecewise bijective and $fi = 1_{H}$; \item $f$  is a
quotient morphism, $N$ is simply connected, and $fi = 1_{H}$.
\end{enumeratein}

\item \label{Ex:9.10:4} Suppose the following diagram of groupoid
morphisms is a pushout $$ \def\labelstyle{\textstyle}
\spreaddiagramcolumns{2pc} %\spreaddiagramrows{2pc}
\diagram {}
\dto_f \rto & {} \dto^g\\ {} \rto & \enddiagram $$  and $f$ is a
quotient morphism.  Prove that $g$ is a quotient morphism.

\item \label{Ex:9.10:4a} Let $f,g: H \to G$ be two groupoid
morphisms. Show how to construct the coequaliser $c:G \to C$  of
$f,g$ as defined in Exercise 6.6.4 of \cite{Brown:1988}. Show how
this gives a construction of the orbit groupoid. [Hint: First
construct the coequaliser $\sigma : \Ob(G) \to Y$ of the functions
$\Ob(f), \Ob(g)$, then construct the groupoid $U_\sigma(G)$, and
finally construct $C$ as a quotient of $U_\sigma(G)$.]

\item \label{Ex:9.10:5} Suppose the groupoid $G$ acts on the
groupoid $\Gamma$ via $w : \Gamma \to \Ob(G)$ as in Exercise
\ref{Ex:9.9:4} of Section 1.  Define the {\em semidirect product
groupoid} $\Gamma \rtimes G$ to have object set $\Ob(\Gamma)$ and
elements the pairs $(\gamma, g) : x \to y$ where $g \in G(wx,wy)$
and $\gamma \in \Gamma (g \cdot x,y)$.  The sum in $\Gamma \rtimes
G$ is given by $(\delta, h) + (\gamma, g) = (\delta + h\cdot
\gamma,hg)$. Prove that this does define a
groupoid, and that the projection % page 9.10-13
$p : \Gamma
\rtimes G \to G, (\gamma, g) \mapsto G$, is a fibration of
groupoids.  Prove that the quotient groupoid $(\Gamma \rtimes
G)/\Ker p$ is isomorphic to $(\pi_{0}\Gamma G)$.

\item \label{Ex:9.10:6} Let $G$ and $\Gamma$ be as in Exercise
\ref{Ex:9.10:5}, and let the groupoid $H$ act on the groupoid
$\Delta$ via $v : \Delta \to \Ob(H)$.  Let $f : G \to H$ and
$\theta : \Gamma \to \Delta$ be morphisms of groupoids such that
$v\theta = \Ob(f)w$ and $\theta (g\cdot \gamma) = (fg)\cdot
(\theta \gamma)$ whenever the left hand side is defined.  Prove
that a morphism of groupoids $(\theta, f)$ is defined by $(\gamma,
g) \mapsto (\theta \gamma, fg)$.  Investigate conditions on $f$
and $\theta $ for $(\theta, f)$ to have the following properties:
\begin{enumeratein} \item injective, \item connected fibres, \item
quotient morphism, \item discrete kernel, \item covering morphism.
\end{enumeratein} In the case that $(\theta, f)$ is a fibration,
investigate the exact sequences of the fibration.

\item \label{Ex:9.10:7} Generalise the Corollary
\ref{Cor:9.10.4:1} from the case of the universal cover to the
case of a regular covering space of $X$ determined by a subgroup
$N$ of $\pi (X,x)$.

\item \label{Ex:9.10:8} Let $1 \to N \to E \to G \to 1$ be an
exact sequence of groups.  Prove that there is an action of $G$ on
a connected groupoid $\Gamma$ and an object $x$ of $\Gamma$ such
that the above exact sequence is isomorphic to the exact sequence
of the fibration $\Gamma \rtimes G \to G$ at the object $x$.
\cite{BrownDanesh-Naruie1}

\item \label{Ex:9.10:9} Let $X^{n}$ be the $n$-fold product of $X$
with itself, and let the symmetric group $S_{n}$ act on $X^{n}$ by
permuting the factors.  The orbit space is called the {\em
$n$-fold symmetric product of $X$} and is written $Q^{n}X$.
 Prove that for $n\geqslant 2$ the fundamental group of $Q^{n}X$
at an image of a diagonal point $(x,\ldots, x)$ is isomorphic to
the fundamental group of $X$ at $x$ made abelian.

\item \label{Ex:9.10:10} Investigate the fundamental groups of
quotients of $X^{n}$ by the actions of various proper subgroups of
the symmetric group for various $n$ and various subgroups. [Try
out first the simplest cases which have not already been done in
order to build up your confidence.  Try and decide whether or not
it is reasonable to expect a general formula.]
\end{enumerate}
\section*{Acknowledgements}
We would like to thank La Monte Yarroll for the major part of the
rendition into Latex of this work. It is hoped to make further
parts  of his efforts available in due course. We would also like
to thank C.D. Wensley for helpful comments.

\end{document}